\documentclass[preprint,12pt]{elsarticle}

\usepackage[T1]{fontenc} 

\usepackage{amsmath} 
\usepackage{amssymb}  
\usepackage{color}
\usepackage{graphicx}
\usepackage{booktabs}
\usepackage{color}
\usepackage{braket} 
\usepackage{algorithm}
\usepackage{cases}
\usepackage{array}

\usepackage{lineno}

\newcommand{\x}{\mathbf{x}}

\newcommand{\boldu}{\mathbf{u}}
\newcommand{\gnom}{\gamma_{\nom}}
\newcommand{\bnom}{\beta_{\nom}}
\newcommand{\sat}{\mathsf{sat}}
\newcommand{\A}{\mathcal{A}}
\newcommand{\U}{\mathcal{U}}

\newcommand{\E}{\mathcal{E}}

\newcommand{\X}{\mathbb{X}}

\newcommand{\dee}{\mathrm{d}}

\newcommand{\dist}{\mathsf{dist}}

\newcommand{\nom}{\mathrm{nom}}
\newcommand{\OCPT}{\mathrm{OCP_T}}
\newcommand{\XM}{\mathbb{X}_{\mathcal{M}}}
\newcommand{\XA}{\mathbb{X}_{\mathcal{A}}}
\newcommand{\M}{\mathcal{M}}

\newcommand{\mbbR}{\mathbb{R}}
\newcommand{\argmin}{\mathrm{argmin}}

\newdefinition{remark}{Remark}
\newtheorem{lemma}{Lemma}
\newtheorem{proposition}{Proposition}
\newtheorem{theorem}{Theorem}
\newproof{problemstatement}{Problem Statement}
\newproof{proof}{Proof}

\journal{...}

\begin{document}

\begin{frontmatter}



\title{MPC without Terminal Ingredients Tailored to the SEIR Compartmental Epidemic Model}

 \author[label1]{Willem Esterhuizen}
 \author[label2]{Philipp Sauerteig} 
  \author[label2]{Stefan Streif}
 \author[label1]{Karl Worthmann}
 
 \affiliation[label1]{organization={Optimization-based Control Group, Institute of Mathematics, Technische Universit\"{a}t Ilmenau},
 	country={Germany} }

 \affiliation[label2]{organization={Automatic Control and System Dynamics, Technische Universität Chemnitz},
 country = {Germany}}

%

\begin{abstract}
We consider the SEIR compartmental epidemic model subject to state and input constraints (a cap on the proportion of infectious individuals and limits on the allowed social distancing and quarantining measures, respectively). 
We present a tailored model predictive control (MPC) scheme without terminal conditions. 
We rigorously show recursive feasibility 
and asymptotic convergence of the MPC closed loop to the continuum of disease-free equilibrium points for suitably designed quadratic running cost and a sufficiently long prediction horizon (forecast window).
Moreover, we establish the viability kernel (a.k.a. the admissible set) as a domain of attraction of the continuum of equilibria. 
\end{abstract}




\end{frontmatter}

\section{Introduction}\label{sec:introduction}

Compartmental models, see for example \cite{Het2000, Brauer2019}, are very popular for describing the spread of infectious diseases. 
In this approach members of a population are assigned to compartments (``susceptible'', ``exposed'', ``infectious'', ``recovered'', ``vaccinated'', etc.) and their movement between them is described by ordinary differential equations. 
The optimal control of these models goes back to at least the 70's. 
To our knowledge, some of the earliest papers on the subject include~\cite{sanders1971quantitative,hethcote1973optimal} --~where optimal epidemic intervention policies are derived via dynamic programming~-- and \cite{gupta1973optimum}, in which an epidemic model is analyzed using 
Pontryagin's maximum principle. 
Since then many papers have appeared on the optimal control of compartmental models, with the SIR (susceptible-infectious-removed) model being by far the most studied, see for example the papers \cite{kirschner1997optimal, Beh00, HanD11, kandhway2014run, agusto2017optimal, BolBS17, godara2021control, borkar2022revisiting, chang2022sparse, BRITTON_20223} (by no means an exhaustive list). 
The recent papers \cite{goreac2022stochastic,freddi2022sir, avram2022optimal} consider the optimal control of stochastic compartmental models, with and without constraints. 

Model predictive control (MPC), see for example \cite{grune2017nonlinear,RawlMayn17}, is an approach where a finite-horizon optimal control problem (OCP) is recursively solved to generate 
a feedback law such that 
the closed-loop trajectory 
is steered to a desired set 
point. 
A well-known strength of MPC is its ability to take 
hard state and input constraints into account. 
In the case of epidemic models, a desired set point might correspond to a disease-free equilibrium, input constraints might model maximal allowed societal interventions 
and state constraints might model a cap on the allowed infection numbers. 

Recently, especially due to the Covid-19 pandemic, many papers have appeared where MPC is applied to 
epidemic models. 
In~\cite{grundel2022much}, the authors derive an age-differentiated compartmental model of Covid-19 and numerically study a number of OCPs concerned with minimising testing and penalising social-distancing measures subject to a hard cap on the hospital capacity. 
They also present a study, where these OCPs are recursively solved in an MPC framework. 
In~\cite{GrunHeyd21b}, the authors extend the model and the analysis proposed in~\cite{grundel2022much} to account for vaccination and study the problem of minimising the effects of social distancing subject to a time-varying constraint on the number of vaccines. 

In~\cite{kohler2020robust}, the authors consider an 8-compartment model of Covid-19 and study an optimal control problem, where the number of fatalities is minimised over two years, without excessive social distancing measures.  
They then study solving this OCP recursively in an MPC framework, and present a robust MPC generalisation to take modelling and measurement uncertainty into account. 
The paper~\cite{selley2015dynamic} applies MPC to a network-based epidemic model, which is of susceptible-infectious-susceptible (SIS) type. 
In this modelling framework individuals are modelled as nodes on a graph that are either susceptible or infectious. 
They control the edges of the graph so as to eradicate the disease while keeping it connected. 
The paper \cite{MorBCetal20} derives an SIR-type model of Covid-19 and studies an MPC problem that aims to reduce the effects of social distancing and the infectious numbers. 
In~\cite{parino2021model} the authors also present a new model of Covid-19 and investigate an MPC concerned with optimal vaccine rollout under constraints on healthcare and social-economic costs. 

All of the so-far mentioned papers that apply MPC to epidemic models are concerned with simulation studies that provide valuable information on how best to manage epidemics. 
However, none of them investigate two very important aspects of an MPC controller: \emph{stability} of the closed-loop system, i.e., the question whether the closed-loop trajectory asymptotically reaches a desired controlled equilibrium point, and \emph{recursive feasibility}. The latter ensures that the finite-horizon OCP is solvable at every time instant provided existence of a solution at the very beginning (initial feasibility). 
To our knowledge, the only papers applying MPC to epidemics that also explicitly look at stability and recursive feasibility are \cite{kohler2018dynamic}, \cite{watkins2019robust}, and our recent paper~\cite{sauerteig2022model}. 

In~\cite{kohler2018dynamic}, the authors model an SIS epidemic as a continuous Markov process and consider the dynamics of the mean-field approximation. 
They present two MPC schemes: one which minimises the disease's spread, the other minimising the cost of containing the disease. 
In the first scheme the recursively solved OCP includes a function, $V_f$, as a terminal cost, whereas the second one uses this function as a terminal constraint.   
Stability and recursive feasibility are guaranteed for both schemes under the assumption that $V_f$ satisfies a Lyapunov-like decrease condition with an initial feasible input.   
The paper~\cite{watkins2019robust} presents robust economic MPC of a stochastic susceptible-exposed-infected-vigilant (SEIV) model, that does not utilise the mean-field approximation. 
As is common in economic MPC, stability of the closed-loop is guaranteed by including a stability constraint function in the finite-horizon OCP. 
The one used in this work ensures that the expectation of the discounted change in the running cost (here, the sum of exposed and infectious nodes) is negative.  

Our previous work~\cite{sauerteig2022model} considers the well-known susceptible-exposed-infectious-removed (SEIR) model with constraints on the state and input and presents an MPC formulation that exploits knowledge of two important sets: the \emph{admissible set} and the \emph{maximal robust positively invariant set} (MRPI), see \cite{DonL13, EstAS20, ester_epidemic_management_2021}. 
We proved that from any state initiating from the admissible set there exists a constraint-satisfying input that asymptotically eliminates the disease, and that this is possible from the MRPI with any input. 
We showed that if one includes a subset of the MRPI as a terminal target set in the finite-horizon OCP, then under some mild assumptions there exists a sufficiently long prediction 
horizon, for which the OCP is recursively feasible and the closed-loop trajectory asymptotically reaches a disease-free equilibrium from \emph{any} initial state located in a suitable subset of the admissible set. 
The inclusion of a target set for the terminal state and a terminal cost in the cost function of the finite-horizon OCP (so-called ``terminal ingredients'') is a well-known way to guarantee stability and recursive feasibility in MPC, see \cite[Ch.~5]{grune2017nonlinear} \cite{mayne2014model}. 

In the current paper, we present an MPC formulation applied to the SEIR model that improves on the work in \cite{sauerteig2022model}. 
We show that terminal ingredients are in fact \emph{not} needed in the finite-horizon OCP to have recursive feasibility and a closed-loop that converges to an equilibrium point.  
Our main result states that a subset of the admissible set, where the infection numbers are not arbitrarily small, is rendered a domain of attraction of the continuum of disease-free equilibrium points under the MPC feedback, provided one uses a suitable quadratic running cost and a sufficiently long prediction horizon. 
Our numerical investigations also suggest that with this formulation the prediction horizon can be much shorter than the one required in \cite{sauerteig2022model}. 
Practically, these facts mean that the iteratively solved finite-horizon OCP can be much simpler (as terminal ingredients are not needed) and the computation time for finding a solution is decreased. 

The outline of the paper is as follows. 
Section~\ref{sec:SEIR_model_and_sets} presents the SEIR model under study and important subsets of the state space that we will consider throughout the paper. 
Section~\ref{sec:MPC} first introduces the finite-horizon OCP and the MPC algorithm we apply to the SEIR model and then presents an in-depth analysis of the OCP. 
It ends with Theorem~\ref{prop:main}, which is the paper's main result. 
Section~\ref{sec:numerics} presents some numerics, and we conclude the paper with Section~\ref{sec:conclusion}. 



\section{Constrained SEIR Model and Sets of Interest}\label{sec:SEIR_model_and_sets}

In this section we present the constrained system under study and introduce a number of subsets of the state space that will be important throughout the paper. 
In particular, we introduce the robustly invariant set, $\XM$, and the admissible set, $\A$. 
Lemma~\ref{lem:facts_of_sets} says that if the state is located in $\XM$ then, regardless of societal intervention, the disease will asymptotically be eradicated and the infection cap will always be respected. 
If the state is located in $\A$, then there exists an intervention such that this is possible. 
These facts regarding the sets will be exploited in the following section where they appear in the MPC formulation. 

\subsection{Constrained SEIR Model}

We consider the following model, see \cite{Het2000},
\begin{subequations}\label{SEIR_IVP}
	\begin{align*}
		\dot{S}(t) 	& = - \beta(t) S(t) I(t), \\
		\dot{E}(t) 	& = \beta(t) S(t) I(t) - \eta E(t),\\
		\dot{I}(t) 	& = \eta E(t) - \gamma(t) I(t),\\
		\dot{R}(t)	& = \gamma(t) I(t), \nonumber \\ 
		(S, E, I, R)(0) & = (S^0, E^0, I^0, R^0)^\top, \nonumber
	\end{align*}
\end{subequations}
$t\geq 0$, where the compartments $S(t)\in\mbbR_{\geq 0}$, $E(t)\in\mbbR_{\geq 0}$, $I(t)\in\mbbR_{\geq 0}$, and $R(t)\in\mbbR_{\geq 0}$ describe the proportions of the population that are \emph{susceptible}, \emph{exposed} (infected, but not yet infectious), \emph{infectious}, or \emph{removed} (recovered, vaccinated, deceased, etc.) at time $t\geq 0$ (measured in days).  
The parameter $\eta^{-1}\in\mathbb{R}_{\geq 0}$ denotes the disease's average \emph{incubation time} in days, whereas $\beta: \mathbb{R}_{\geq 0} \to \mathbb{R}_{> 0}$ and $\gamma: \mathbb{R}_{\geq 0} \to \mathbb{R}_{> 0}$ are time-varying manipulatable inputs representing \emph{social distancing and hygiene concepts}, and \emph{quarantine measures}, respectively.
We assume that these inputs are subject to constraints, $\beta(t) \in [\beta_{\min}, \beta_{\nom}]$, $\gamma(t) \in [\gamma_{\nom}, \gamma_{\max}]$, for all $t\geq 0$, where $0<\beta_{\min} < \beta_{\nom}$ and $0 < \gamma_{\nom} < \gamma_{\max}$. 
The values $\beta_{\nom}$ and $\gamma_{\nom}$ denote the ``nominal values'' of these inputs, (the nominal \emph{contact rate} and \emph{recovery/removal rate} of the disease, respectively) in the absence of societal intervention, whereas
$\beta_{\min}$ and $\gamma_{\max}$ are the limits of allowed social distancing and quarantining measures, respectively. 

\sloppy
We aim to maintain a hard infection cap via the state constraint,
$I(t) \leq I_{\max}$ for all $t\geq 0$, where $0 < I_{\max} \leq 1$ represents the \emph{maximal hospital capacity}. 
Without loss of generality, we assume that the initial state $(S^0,E^0,I^0,R^0)^\top~\in~[0,1]^4\subset\mathbb{R}^4$ satisfies $S^0 + E^0 + I^0 + R^0 = 1$. 
Noting that $\dot S + \dot E + \dot I + \dot R = 0$  we see that $S(t) + E(t) + I(t) + R(t) = 1$ for all $t\geq 0$ for any pair of inputs $\beta$, $\gamma$. 
Because the compartment $R$ has no influence on the other compartments we will ignore it and only consider the three-dimensional system involving $\dot S$, $\dot E$ and $\dot I$ in the sequel.

\fussy
To ease our notation we write the constrained system as follows,
\begin{numcases}{\mathcal{S}:}
	& $\dot{\x}(t)  = f(\x(t),\boldu(t))$, \nonumber \\
	& $\x(0) = \x^0$, \nonumber\\
	& $\x(t) \in \mathbb{X}$,\nonumber\\
	& $\boldu(t) \in \mathbb{U}$, \nonumber
\end{numcases}
$t\geq 0$, where $\x:=(x_1,x_2,x_3)^\top\in\mbbR^3$ denotes the state, with $x_1=S$, $x_2=E$, $x_3=I$, and $\boldu:=(u_1,u_2)^\top\in\mathbb{R}^2$ denotes the control, with $u_1 = \beta$, $u_2 = \gamma$. 
The system equations $f:\mathbb{R}^3\times\mathbb{R}^2\rightarrow \mathbb{R}^3$ read, 
\begin{equation}
	f(\x,\boldu) = \left(
	\begin{array}{c}
		-u_1x_1x_3\\ 
		u_1x_1x_3 - \eta x_2\\ 
		\eta x_2 - u_2 x_3
	\end{array}
	\right), \label{eq:SIER_f}
\end{equation}
and the initial state is specified by $\x^0 := (x_1^0,x_2^0,x_3^0)^\top = (S^0,E^0,I^0)^\top~\in~[0,1]^3$. 
We define the function $g:\mathbb{R}^3\rightarrow\mathbb{R}$ as follows,
\[
g(\x):= x_3 - I_{\max}, 
\]
and let the set $\mathbb{X}$ read,
\begin{equation}
	\X := \{\x\in[0,1]^3: \left( g(\x) \leq 0 \right) \wedge \left( x_1 + x_2 + x_3 \leq 1 \right)\}. \nonumber
\end{equation}
Finally, $\mathbb{U}$ is defined as follows,
\begin{equation}
	\mathbb{U} := [\beta_{\min}, \beta_{\nom}]\times[\gamma_{\nom}, \gamma_{\max}]\subset\mathbb{R}_{> 0}\times\mbbR_{> 0}.\nonumber
\end{equation}
We let 
\[
\mathcal{U}_{\infty} = \mathcal{L}([0,\infty),\mathbb{U}),
\] 
be the set of Lebesgue measurable functions mapping $[0,\infty)$ to $\mathbb{U}$. 
Moreover, with $T\in\mbbR_{\geq 0}$, we denote the finite horizon counterpart by,
\[
	\mathcal{U}_T = \mathcal{L}([0,T],\mathbb{U}). 
\]
Occasionally we will slightly abuse our notation by using $\boldu$ to indicate an element of a function space ($\U_{\infty}$ or $\U_{T}$), or the space $\mbbR^2$, but the context should cause no confusion. 

The function $f$ is smooth with respect to $\x$ and Lebesgue measurable with respect to $t$. 
Moreover, for any initial state $\x^0=(x_1^0,x_2^0,x_3^0)^\top$ satisfying $x_1^0 + x_2^0 + x_3^0 \leq 1$, the solution of the system remains in the compact set $[0,1]^3$ for all $t\geq 0$. 
Therefore, for every such $\x^0$ and every $\boldu\in\mathcal{U}_{\infty}$ there exists a unique locally absolutely continuous curve, $\x:\mathbb{R}_{\geq 0}\rightarrow\mathbb{R}^3$, that satisfies $\dot{\x}(t) = f(\x(t),\boldu(t))$ for a.e. $t\in\mathbb{R}_{\geq 0}$. 
See for example \cite[App. A]{Sontag2013}. 
By $\x(t;\x^0,\boldu) := (x_1(t;\x^0,\boldu),x_2(t;\x^0,\boldu), x_3(t;\x^0,\boldu))^\top$ we denote this solution at time $t\in\mbbR_{\geq 0}$, with input $\boldu\in\U_{\infty}$, initiating at $t=0$ from $\x^0\in[0,1]^3$. 

\subsection{The admissible set and disease-free equilibria}

\begin{figure}[htb]
	\centering
	\includegraphics[scale=1]{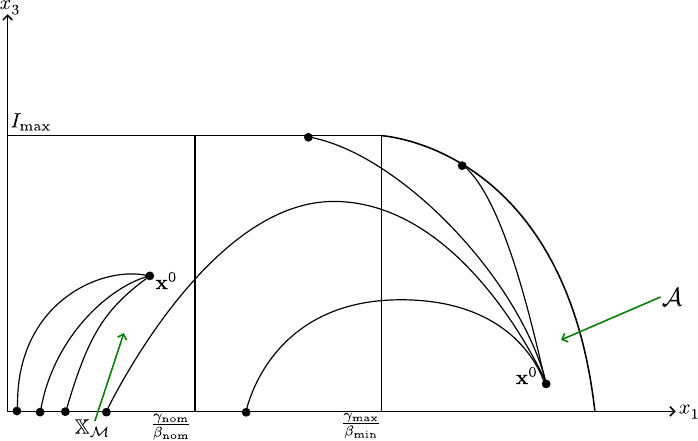}
	\caption{Admissible set, $\A$, and invariant box, $\XM$. 
    From $\x^0\in\A$, with an admissible input, it is possible for the state to always satisfy the infection cap and to asymptotically settle at a point where $x_1\leq \frac{\gamma_{\max}}{\beta_{\min}}$. 
    From $\x^0\in\XM$ the state will satisfy the infection cap for any admissible input and settle at a point where $x_1\leq \frac{\gamma_{\nom}}{\beta_{\nom}}$.}
	\label{fig:Lemma_explanation_1}
\end{figure}

Throughout the paper we will consider the following box, 
\begin{align}
    \mathbb{X}_{\M} := \{\x\in\mathbb{X} : x_1 \leq \overline{x}_1, x_2\leq \overline{x}_2\}, \label{eq:def_set_XM}
\end{align}
where,
\[
	\overline{x}_1 := \frac{\gamma_{\nom}}{\beta_{\nom}},\quad  \overline{x}_2 := \frac{\gamma_{\nom}I_{\max}}{\eta},
\]
as well as the \emph{admissible set}, \cite{DonL13}, (also known as the viability kernel, \cite{aubin2011viability}),
\begin{align}
	\A := \{\x^0\in \X: \exists \boldu\in\U_{\infty},\,\,\x(t;\x^0,\boldu)\in \X\,\,\forall t\geq 0\}. \label{eq:def_set_A}
\end{align}
The following lemma summarizes some important facts about these sets. 
It says that from the admissible set it is possible, via an input that satisfies the input constraints, to asymptotically eradicate the disease while always satisfying the state constraint. 
From the box $\XM$ this is possible regardless of the input (thus, $\XM$ is a robust invariant set contained in $\X$, see for example \cite{EstAS20}). 
In other words: The set~$\XM$ is a safe set, meaning that even without interventions the hard infection cap will be maintained. 
For this reason, for a given initial value $\x^0 \in \mathcal{A}$ we would like the state to end up in the safe set~$\XM$.
Refer to Figure~\ref{fig:Lemma_explanation_1} for further clarity on Lemma~\ref{lem:facts_of_sets}. 
\begin{lemma}\label{lem:facts_of_sets}
    The following assertions hold:
    \begin{enumerate}
	\item For every $\x^0 \in \A$, defined in~\eqref{eq:def_set_A}, there exists $\boldu \in \U_{\infty}$ such that $\x(t;\x^0,\boldu) \in \X$ holds for all $t\geq 0$ and 
            \begin{align*} 
                & \lim_{t \to \infty} x_1(t; \x^0,\boldu) \leq \frac{\gamma_{\max}}{\beta_{\min}} \in[0,1], \\ 
                & \lim_{t \to \infty}x_2(t; \x^0,\boldu) = \lim_{t \to \infty} x_3(t; \x^0,\boldu) = 0.
            \end{align*}
		\item For any $\x^0\in\XM$, defined in~\eqref{eq:def_set_XM} and for all $\boldu\in\U_{\infty}$, the state stays within this set, i.e., $\x(t;\x^0,\boldu)\in\XM \subseteq \X$ for all $t \geq 0$, i.e., $\XM \subseteq \A$. 
            Moreover, the following holds: 
    		\begin{align*} 
    			& \lim_{t \to \infty} x_1(t; \x^0,\boldu) \leq \frac{\gamma_{\nom}}{\beta_{\nom}} \in \left[ 0, \frac{\gamma_{\max}}{\beta_{\min}} \right], \\
    			& \lim_{t \to \infty} x_2(t; \x^0,\boldu) = \lim_{t \to \infty} x_3(t; \x^0,\boldu) = 0.
    		\end{align*}
	\end{enumerate}
\end{lemma}
\begin{proof}
	Detailed proofs of~1. and~2. have been given in \cite[Lem. 1]{sauerteig2022model} and \cite[Lem. 2]{sauerteig2022model}, respectively. 
	For sake of self-containedness we sketch the main ideas here. 
	The only thing not explicitly stated in~\cite{sauerteig2022model} is the terminal bound on $\lim_{t \to \infty} x_1(t; \x^0, \boldu)$ both in~1. and~2. 
	\begin{enumerate}
		\item In a first step, we showed that the limits $\lim_{t \to \infty} x_i(t; \x^0, \boldu) \in [0,1]$ exist for $i \in \{1,2,3\}$ and $\lim_{t \to \infty} x_2(t; \x^0, \boldu) = \lim_{t \to \infty} x_3(t; \x^0, \boldu) = 0$ for all feasible $\x^0$ and $\boldu$. 
		      Based on this, one can deduce that for all $\x^0 \in \mathcal{A}$ there exists a $\boldu \in \mathcal{U}_\infty$ such that $\x(t; \x^0, \boldu) \in \mathbb{X}$ for all $t \geq 0$ and $\lim_{t \to \infty} x_2(t; \x^0, \boldu) = \lim_{t \to \infty} x_3(t; \x^0, \boldu) = 0$. 
    		Moreover, due to 
    		\begin{align*}
    			\frac{\mathrm{d}}{\mathrm{d} t} \left( x_2(t) + x_3(t) \right) 
    			\begin{cases}
    				> 0 \quad & \text{if } x_1(t) > \gamma(t) / \beta(t), \\
    				= 0 \quad & \text{if } x_1(t) = \gamma(t) / \beta(t), \\ 
    				< 0 \quad & \text{if } x_1(t) < \gamma(t) / \beta(t), 
    			\end{cases} 
    		\end{align*}
    		the terminal value $\lim_{t \to \infty} x_1(t; \x^0, \boldu)$ has to satisfy
    		\begin{align*}
    			\lim_{t \to \infty} x_1(t; \x^0, \boldu) \leq \lim_{t \to \infty} \frac{\gamma(t)}{\beta(t)} \leq \frac{\gamma_{\max}}{\beta_{\min}} 
    		\end{align*}
    		in order for $\lim_{t \to \infty} x_2(t; \x^0, \boldu) = \lim_{t \to \infty} x_3(t; \x^0, \boldu) = 0$ to hold. 
		\item In \cite[Lem.~2]{sauerteig2022model}, we showed that for $\x \in \XM$, the states~$x_2$ and~$x_3$ decay exponentially under any feasible control input by showing that these states are bounded from above by a linear system with exponential decay. 
            Thus, the set $\XM$ is forward invariant under any feasible control input. 
            The remainder follows from the definition of~$\XM$ and the exponential decay of~$x_2$ and~$x_3$. 
	\end{enumerate}
\end{proof}
We let $\E:= \{\x\in\X: x_1\in[0,1], x_2 = 0, x_3 = 0\}$ denote the continuum of disease-free equilibrium points, and let
\[
\E_{\nom} = \Set{\x\in\E: x_1\leq \overline{x}_1},
\]
denote those equilibrium points that are contained in $\XM$. 
Finally, in the sequel we will consider the following subset of $\A$,
\begin{equation}
\A^{\prime} :=  \A\setminus \mathcal{N}_{\varepsilon}, \label{eq:def_set_A_prime}
\end{equation}
where,
\[
	\mathcal{N}_{\varepsilon} :=\{\x^0\in\A : \left( x_1^0 \geq \bar{x}_1 \right)\, \land  \left( x_3^0 \leq \varepsilon \right) \},
\]
with $0 < \varepsilon \leq 1$ small. 
Thus, $\A^{\prime}$ excludes points outside $\XM$ for which the initial number of infectious individuals are arbitrarily small. 
As will be clarified, this is to ensure that the state can reach the invariant box $\XM$ in finite time, a fact we need to prove the uniform boundedness of the infinite horizon value function (Proposition~\ref{lem:V_inf_uni_bounded_on_A}).

Considering $\A^{\prime}$ instead of $\A$ in the sequel rigorously addresses a peculiarity that arises from modelling a finite number of interacting humans with a continuous state. 
In plain language, $\A^{\prime}$ excludes meaningless cases where the initial infected population consists of a fraction of a person. 
However, if the initial condition satisfies $\x^0\in$  $\mathcal N_{\varepsilon}$, then there exists a finite time $t_0(\x^0) > 0$ such that $\x(t_0(\x^0); x^0, \boldu_{\nom}) \in A^{\prime}$ with $\x(t; x^0, \boldu_{\nom}) \in \X$ for all $t \in [0, t_0(x_0)]$. 
In other words, those academic cases require no special treatment.

\section{MPC without Terminal Ingredients}\label{sec:MPC}

In this section we first introduce the recursively solved finite-horizon optimal control problem and then specify the MPC scheme we apply to the constrained SEIR model, $\mathcal{S}$, in Algorithm~\ref{alg:mpc}. 
We then present an in-dpeth analysis of the OCP and MPC scheme, arriving at our main result, Theorem~\ref{prop:main}. 

\subsection{The Finite-Horizon Optimal Control Problem}

Consider the following optimal control problem, where $T\in\mbbR_{\geq 0}\cup\{+\infty\}$:
\[
\OCPT : \quad \min_{\boldu \in \U_T(\x^0)} J_T(\x^0,\boldu),
\]
with cost functional, $J_T:\mbbR^3\times \U_T \rightarrow \mbbR_{\geq 0}$, given by
\[
	J_T(\x^0,\boldu) := \int_0^T \ell(\x(s;\x^0,\boldu),\boldu(s))\,\dee s, 
\]
and where $\U_T(\x^0)$ (resp. $\U_{\infty}(\x^0)$) denotes all control functions $\boldu \in \U_T$ (resp. $\boldu \in\U_{\infty}$) for which $\x(t;\x^0,\boldu)\in \X$ for all $t \in [0,T]$ (resp. $t \in [0,\infty)$). 
The stage cost is given by $\ell : \mathbb{R}^3 \times \mathbb{R}^2 \to \mathbb{R}_{\geq 0}$ and we will consider the following one throughout the paper, with $\lambda\in(0,1]$,
\begin{equation}
	\ell(\x,\boldu) := \lambda\Vert \x \Vert_Q^2 + (1 - \lambda)\Vert \boldu - \boldu_{\nom}\Vert^2, \label{quadratic_ell}
\end{equation}
where $\boldu_{\nom} := (\beta_{\nom},\gamma_{\nom})^\top$ and 
$Q = \left(
\begin{array}{ccc}
	0 & 0 & 0\\
	0 & 1 & 0\\
	0 & 0 & 1
\end{array}
\right)$. 
Here, $\Vert \x \Vert_Q^2 := \x^\top Q \x$, 
Moreover, we define,
\[
\ell^{\star}(\x) := \min_{\boldu\in\mathbb{U}}\ell(\x,\boldu) = \lambda\Vert \x \Vert_Q^2.
 \]
Thus, the cost functional penalises deviation of the state from the continuum of disease-free equilibria as well as deviation of the control from its ``nominal'' value (i.e., it minimises societal intervention).
The parameter $\lambda\in(0,1]$ may be used to weight the importance of these two objectives. 
Note that with this choice of $\ell$ we do not aim for a specific equilibrium point in $\E$ (because $x_1$ is not penalised in the cost) and that $\ell^{\star}$ is \emph{not} positive definite.  
However, as Theorem~\ref{prop:main} in the next subsection will show, this presents no issues in deducing recursive feasibility of $\OCPT$ and an MPC closed-loop that asymptotically approaches the continuum of disease-free equilibrium points. 
For more results on MPC with indefinite stage costs we refer the interested reader to the papers \cite{berberich2018indefinite, kohler2023stability}.   

We indicate the value function $V_T:\mathbb{R}^3\rightarrow \mathbb{R}_{\geq 0}\cup \{+\infty\}$, $T\in\mbbR_{\geq 0}\cup\{+\infty\}$ by,
\[
	V_T(\x^0):= \min_{\boldu\in\U_T(\x^0)} J_T(\x^0,\boldu), 
\]
and, with $C\in\mathbb{R}_{\geq 0}$, the sublevel sets by,
\[
	V_T^{-1}[0,C] := \{\x\in\mathbb{R}^3 : V_T(\x) \leq C\}. 
\]
Throughout the paper we assume that if $\mathcal{U}_T(\x^0) \neq \emptyset$ for an arbitrary $\x^0~\in~\mbbR^3$ and $T\in\mbbR_{\geq 0}\cup\{+\infty\}$, then there exists an optimal control (a global minimiser) in  $\mathcal{U}_T(\x^0)$. 
For existence results in optimal control we refer the reader to \cite[Ch.9-Ch.16]{cesari2012optimization} and \cite[Ch.4]{Macki2012}. 
Furthermore, if $\mathcal{U}_T(\x^0) = \emptyset$ for an arbitrary $\x^0$ and $T$, we take the convention that $V_T(\x^0) = +\infty$. 

\subsection{The MPC scheme without terminal ingredients}

The MPC algorithm we utilise in this paper is given in Algorithm~\ref{alg:mpc}. 
Note the absence of a terminal cost and terminal constraint set. 
Thus, starting at an initial state $\x^0\in\A^{\prime}$, the finite horizon problem $\OCPT$ is solved, the first part for $t\in[0,\delta]$ is implemented, the state is updated, and $\OCPT$ is solved again. 
This iteration is done for all time, producing the so-called \emph{MPC feedback}, $\mu_{T,\delta} : [0,\delta)\times \mathbb{X} \rightarrow \mbbR^2$, given by $\mu_{\delta, T}(t,\x^0) = \boldu^{\star}(t)$, for $t\in[0,\delta]$.  
\begin{algorithm}[h]
	\caption{Continuous-time MPC without terminal ingredients}
	\noindent \textbf{Input}: control horizon, $\delta \in \mathbb{R}_{>0}$, number of steps, $N \in \mathbb{N}$, initial state, $\x^0 \in \mathcal{A}^{\prime}$\\
	\textbf{Set}: prediction horizon, $T \leftarrow N \delta$
	\begin{enumerate}
		\item Find a minimizer $\boldu^\star \in \argmin_{\boldu \in \mathcal{U}_T(\x^0)} J_T(\x^0,\boldu)$.
		\item Implement $\boldu^\star(\boldsymbol \cdot)$, for $t \in [0, \delta]$.
		\item Set $\x^0\leftarrow \x(\delta;\x^0,\boldu^{\star})$ and go to step 1).
	\end{enumerate}
	\label{alg:mpc}
\end{algorithm}

\subsection{Recursive Feasibility and Convergence to Equilibrium Points}

In this section we first establish a number of facts concerning $\OCPT$. 
Then we state the paper's main result, Theorem~\ref{prop:main}. 
It states that for any $\delta >0$ there exists a $T \geq \delta$ such that $\OCPT$ is recursively feasible in the MPC loop of Algorithm~\ref{alg:mpc}, and such that the closed loop asymptotically approaches an equilibrium point on $\E_{\nom}$. 

The following lemma states that, with the nominal control (i.e., with no societal intervention), the infinite horizon cost functional is uniformly bounded on the invariant box, $\XM$. 
\begin{lemma}\label{lem:J_inf_uni_bounded}
	There exists a $C\in\mathbb{R}_{\geq 0}$ such that,
	\[
	J_{\infty}(\x^0,\boldu_{\nom})\leq C,\,\,\forall \x^0\in\XM,
	\]
    $\XM$ defined in \eqref{eq:def_set_XM}. 
	Here, $\boldu_{\nom}\in\U_{\infty}$ indicates the constant function $\boldu_{\nom}(t) = (\beta_{\nom},\gamma_{\nom})^\top$ for all $t\geq 0$. 
\end{lemma}
\begin{proof}
	The proof follows the same arguments as that of \cite[Lem. 3]{sauerteig2022model}, with only a slight generalisation to accommodate the weighting parameter $\lambda$. 
	See the appendix for the proof. 
\end{proof}
The following proposition states that \emph{cost controllability}, see for example \cite{coron2020model,grune2010analysis, kohler2023stability} and \cite[Ch.~6]{grune2017nonlinear}, holds over the invariant box, $\XM$. 
This is a key property of the system that we will exploit in the proof of Theorem~\ref{prop:main}. 
\begin{proposition}\label{prop:cost_control}
	There exists a $\rho\in\mathbb{R}_{\geq 0}$ such that,
	\[
	V_{\infty}(\x^0)\leq \rho \ell^{\star}(\x^0),\,\,\forall \x^0\in\XM,
	\]
 $\XM$ defined in \eqref{eq:def_set_XM}. 
\end{proposition}
\begin{proof}
	\fussy	
	To ease our notation let $\x(t) = (x_1(t), x_2(t), x_3(t))^\top = \x(t;\x^0,\boldu_{\nom})$ denote the ``nominal'' solution initiating from $\x^0 = (x_1^0, x_2^0, x_3^0)^\top\in\XM$. 	
	From \cite[Lem.~2]{sauerteig2022model} we have that, for any $\x^0\in\XM$,
	\begin{align*}
		\left(
			\begin{array}{c}
				\dot x_2(t)\\
				\dot x_3(t)				
			\end{array}
		\right)
		\leq 
		A
		\left(
		\begin{array}{c}
			x_2(t)\\
			x_3(t)				
		\end{array}
		\right),
		\quad \forall\,  t\geq 0, \quad (x_2(0),x_3(0))^\top = (x_2^0,x_3^0)^\top,
	\end{align*}
	where
	$A$ is a Hurwitz matrix. 
	Therefore, there exist constants, $\Gamma>0$ and $\eta >0$ such that,
	\begin{align*}
		\left\Vert
		\left(
		\begin{array}{c}
			x_2(t)\\
			x_3(t)				
		\end{array}
		\right)
		\right\Vert
		\leq 
		\Gamma e^{-\eta t}
		\left\Vert
		\left(
		\begin{array}{c}
			x_2(0)\\
			x_3(0)				
		\end{array}
		\right)
		\right\Vert
		\quad \forall\,  t\geq 0,\quad (x_2(0),x_3(0))^\top = (x_2^0, x_3^0)^\top,
	\end{align*}
	$\| \cdot \|$ indicating the Euclidean norm. 
	Therefore,
	\begin{align*}
		x_2^2(t) + x_3^2(t) \leq
			\left( \Gamma e^{-\eta t}\left\Vert
			\left(
			\begin{array}{c}
					x_2(0)\\
					x_3(0)				
				\end{array}
			\right)
			\right\Vert
			\right)^2
			= \left( \Gamma e^{-\eta t} \right)^2
			\left(
			(x_2^0)^2 + (x_3^0)^2
			\right).
		\end{align*}
	Therefore, for all $\x^0\in\XM$,
	\begin{align*}
		V_{\infty}(\x^0) \leq J_{\infty}(\x^0, \boldu_{\nom}) & = \lambda \int_0^{\infty} 	\left( x_2^2(s) + x_3^2(s)\right)\,\, \dee s\\
		& \leq \lambda \Gamma^2\left(
		(x_2^0)^2 + (x_3^0)^2
		\right)
		\int_0^{\infty} e^{-2\eta s}  \,\, \dee s \\
		& = \frac{\lambda\Gamma^2}{2\eta}\left(
		(x_2^0)^2 + (x_3^0)^2
		\right)
		= \rho \ell^{\star}(\x^0),
	\end{align*}
	where $\rho := \frac{\lambda\Gamma^2}{2\eta}$.
		\fussy
\end{proof}
\begin{remark}
	The statement in Proposition~\ref{prop:cost_control} does not hold for initial states with $x_1 > \bar{x}_1$. To see this, consider an equilibrium point $\x^{\mathrm{eq}} = (x^{\mathrm{eq}}_1, 0, 0)^\top$ with $x^{\mathrm{eq}}_1 > \bar{x}_1$ and $\varepsilon < x^{\mathrm{eq}}_1 - \bar{x}_1$. %
	Then, for all $\x^0  \in \mathcal{N}_\varepsilon(\x^{\mathrm{eq}})$ with $\min \{x_2,x_3\} > 0$, the value $V_{\infty}(\x^0)$ is positive since $x_1 - \bar{x}_1$ many people have to get infected before we can switch to $\boldu_{\nom}$. %
	However, $\ell^\star(\x^0)$ tends towards 0 for $x_2^0,x_3^0 \to 0$, and so $\rho$ would tend to infinity. 
\end{remark}
\fussy
We now show that $V_{\infty}$ is uniformly bounded on the subset of the admissible set, $\A^{\prime}$. 
\begin{proposition}\label{lem:V_inf_uni_bounded_on_A}
    There exists a $C^\prime\in\mathbb{R}_{\geq 0}$ such that,
	\[
	V_{\infty}(\x^0) \leq C^\prime,\,\forall\,  \x^0\in\A^{\prime},
	\]
    $\A^\prime$ defined in \eqref{eq:def_set_A_prime}.
In other words, there exists a $C^{\prime}\in\mathbb{R}_{\geq 0}$ such that $\A^{\prime} \subseteq V_{\infty}^{-1}[0,C^{\prime}]\subseteq V_{T}^{-1}[0,C^{\prime}]$, $T < \infty$. 
\end{proposition}
\begin{proof}
	\sloppy
	The idea of the proof is sketched in Figure~\ref{fig:Lemma_explanation}. 
	First, because the state is located in $\A$, we utilise an input $\tilde{\boldu}$ with which the resulting state trajectory satisfies the infection cap for all time until $x_1 \leq \overline x_1$ and $\dot x_3 = 0$. 
	We then switch to a special feedback that keeps $x_3$ constant until the state enters the invariant box, $\XM$, at which point we switch to the nominal input, $\boldu_{\nom}$, which drives the state to an equilibrium point. 
	We argue that if $\x^0\in\A^{\prime}$ then $\XM$ is reached in finite time.
	\begin{figure}[htb]
		\centering
		\includegraphics[scale=1]{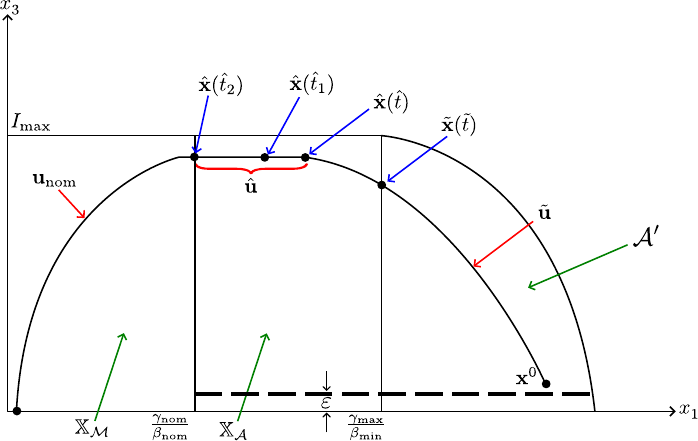}
		\caption{Sketch of the proof of Proposition~\ref{lem:V_inf_uni_bounded_on_A}.}
		\label{fig:Lemma_explanation}
	\end{figure}
	
	\noindent
	Without loss of generality, consider an initial state $\x^0\in\mathcal{A}^{\prime}$ with $x_1 > \underline{x}_1 = \frac{\gamma_{\max}}{\beta_{\min}}$. 
	Recall from Lemma~\ref{lem:facts_of_sets} that there exists a $\tilde \boldu \in \mathcal{U}_{\infty}$ such that $g(\x(t;\x^0,\tilde\boldu)) \leq 0$ for all $t\geq 0$ and $x_1^{\infty}(\x^0,\tilde\boldu) \leq \frac{\gamma_{\max}}{\beta_{\min}} \in[0,1]$, $x_2^{\infty}(\x^0,\tilde\boldu) = x_3^{\infty}(\x^0,\tilde\boldu) = 0$. 
	To ease our notation, let $\tilde \boldu(t) := (\tilde u_1(t), \tilde u_2(t))^\top$ and $\tilde \x(t) := (\tilde x_1(t), \tilde x_2(t), \tilde x_3(t))^\top$, where $\tilde x_i(t) := x_i(t;x,\tilde \boldu)$, $i=1,2,3$. 
	Because $\tilde x_1$ is strictly decreasing, there exists a $\tilde t \in\mbbR_{>0}$ for which $\tilde x_1(\tilde t) = \underline{x}_1 > 0$. 
	If $x_1 > \underline{x}_1$ then $\min_{\boldu\in\mathbb{U}}\dot x_2 + \dot x_3 = \left( \beta_{\min} x_1 - \gamma_{\max} \right)x_3 > \left( \gamma_{\max} - \gamma_{\max} \right)x_3 = 0$, thus $\tilde x_2(t) + \tilde x_3(t)$ is strictly increasing over $[0,\tilde t]$. 
	Define the constant,
	\[
		K := \min_{t\in[0,\tilde t]} \frac{\tilde x_3(t)}{\tilde x_2(t) + \tilde x_3(t)}, 
	\]
	from where we deduce that,
	\begin{equation}
		K(\tilde x_2(t) + \tilde x_3(t))\leq \tilde x_3(t), \quad \forall t\in[0,\tilde t]. \label{eq:K_equation}
	\end{equation}
	Because $\dot x_3 \geq -\gamma_{\max}x_3$, we have $\tilde x_3(t) \geq e^{-\gamma_{\max}t}x_3^0$ for all $t\in[0,\tilde t]$. 
	Recall that $\tilde x_2(t) + \tilde x_3(t) \leq 1 - \tilde x_1(t) \leq 1 - \underline x_1$, $t\in[0,\tilde t]$.  
	Thus, 
	\[
		K  \geq \frac{e^{-\gamma_{\max} \tilde t}x_3^0}{1 - \underline x_1}.
	\]
	Using Equation~\eqref{eq:K_equation}, and the fact that $\tilde x_2(t) + \tilde x_3(t)$ is strictly increasing, we see that $\dot{ \tilde x}_1(t) = -\tilde u_1(t) \tilde x_1(t) \tilde x_3(t) \leq -\beta_{\min} K(\tilde x_2(t) + \tilde x_3(t))\tilde x_1(t) \leq -\beta_{\min} K(x_2^0 + x_3^0)\tilde x_1(t) \leq -\tilde{K}\tilde x_1(t)$, where we have defined,
	\[
		\tilde K := \beta_{\min}\frac{e^{-\gamma_{\max} \tilde t}x_3^0}{1 - \underline x_1}(x_2^0 + x_3^0).
	\]
	Because $\x^0\in\A^{\prime}$ we have $x_3^0(x_2^0 + x_3^0) \geq \varepsilon_2 > 0$, with $\varepsilon_2 < \varepsilon$, thus $\tilde K  \geq \beta_{\min} \frac{e^{-\gamma_{\max} \tilde t}}{1 - \underline x_1} \varepsilon_2 > 0$. 
	Thus, from $\tilde x_1(t) \leq e^{-\tilde K t} x_1^0$, we see that $\tilde t$ is finite. 
	
	Because $x_3^{\infty}(\x^0,\tilde \boldu) = 0$ there exists an interval $[\tilde t, \hat t]$, with $\hat t > \tilde t$, over which $\tilde x_3$ is strictly increasing, and $\dot{\tilde x}_3(\hat t) = 0$. 
	The time $\hat t$ must be finite, for otherwise $x_3^{\infty}(\x^0,\tilde \boldu) \neq 0$, which is a contradiction.	
	At $t = \hat t$, $\dot{\tilde x}_3(\hat t) = 0$, which implies $\eta \tilde x_2(\hat t) - \tilde u_2(\hat t) \tilde x_3(\hat t) = 0$. 
	Thus, $\tilde x_2(\hat t) = \frac{\tilde u_2(\hat t) \tilde x_3(\hat t)}{\eta} \leq \frac{\gamma_{\max} I_{\max}}{\eta}$. 
	Thus, at $t = \hat t$ the state has reached the box $\XA :=  [0,\frac{\gamma_{\max}}{\beta_{\min}}]\times [0,\frac{\gamma_{\max}I_{\max}}{\eta}] \times [0, I_{\max}]$. 
	Now, for all $t\in[\hat t,\infty)$ switch to the feedback $\hat \boldu(\x(t)) := (\hat u_1(\x(t)), \hat u_2( \x ( t ) ) )^\top$,
	\[
		\hat u_1(\x(t)) := \sat\left(\frac{\eta x_2(t)}{x_1(t)x_3(t)},\beta_{\min},\bnom \right), \hat u_2(\x(t)) := \sat\left(\frac{\eta x_2(t)}{x_3(t)},\gnom, \gamma_{\max} \right)
	\]
	To ease our notation we will drop $\hat\boldu$'s dependence on the state and let $\hat \x(t) := (\hat x_1(t), \hat x_2(t), \hat x_3(t))^\top$, where $\hat x_i (t) := x_i(t;\tilde \x(\hat t), \hat \boldu)$, $t\in[\hat t,\infty)$, $i=1,2,3$, denotes the closed loop solution with feedback $\hat \boldu$, from $\tilde\x(\hat t)$. 
	Thus, $\tilde \x(\hat t) = \hat \x(\hat t)$. 
	Recall that $\dot{\hat x}_3(\hat t) = \dot{\tilde x}_3(\hat t) = 0$, so there exists a $u_2\in[\gamma_{\nom}, \gamma_{\max}]$ such that $\dot{x}_3(\hat t) = 0$ and we can thus deduce that $\frac{\eta x_2(\hat t)}{x_3(\hat t)} \in [\gamma_{\nom}, \gamma_{\max}]$.  
	Moreover, because $\dot{x}_3(\hat t) = 0$, $ \min_{(u_1,u_2)\in\mathbb{U}}\dot {x}_2(\hat t) = \min_{(u_1,u_2)\in\mathbb{U}}\dot {x}_2(\hat t) +  \dot{x}_3(\hat t) = (\beta_{\min} \hat x_1(\hat t)-\gamma_{\max})\hat x_3(\hat t) \leq 0$ because $\hat x_1(\hat t) \leq \frac{\gamma_{\max}}{\beta_{\min}}$.
	Thus, $\frac{\eta x_2(\hat t)}{x_1(\hat t)x_3(\hat t)} \in [\beta_{\min}, \beta_{\nom}]$.
	
	Now, there exists a time interval, $[\hat t, \hat t_1]$, over which $\hat x_2$ and $\hat x_3$ remain constant while $\hat x_1$ is still decreasing, and therefore $\hat u_1$ increases and $\hat u_2$ remains constant. 
	Let $\hat t_1$ be the time for which $\hat u_1(\hat t) = \beta_{\nom}$. 
	We see that $\hat x_1(\hat t_1) = \frac{\eta \hat x_2(\hat t)}{\beta_{\nom}\hat x_3(\hat t)} > 0$. 
	Because $\hat x_3$ remains constant on this interval, and because $x_3$ was strictly increasing before $\hat t$, we have the estimates,
	\[
		\hat x_1(t) \leq e^{-\beta_{\min} \hat x_3(\hat t)(t - \hat t)}\hat x_1(\hat t),\quad t\in[\hat t, \hat t_1],
	\]
	and
	\[
	\hat x_3(\hat t) \geq e^{-\gamma_{\max} \hat t}x_3^0. 
	\]	
	Moreover, again using the fact that $\x^0\in\A^{\prime}$, we have $x_3^0(x_2^0 + x_3^0) \geq \varepsilon_2$, and therefore, $x_3^0 \geq \frac{\sqrt{\varepsilon_2} - 1}{2} >0$.   
	Therefore, $\beta_{\min} \hat x_3(\hat t) \geq \beta_{\min}e^{-\gamma_{\max} \hat t}x_3^0 \geq \beta_{\min}e^{-\gamma_{\max} \hat t}\frac{\sqrt{\varepsilon_2} - 1}{2}> 0$. 
	Thus, $\hat t_1$ is finite. 

	Next, there exists an interval, $[\hat t_1, \hat t_2]$, over which $\hat u_1 \equiv \beta_{\nom}$, and thus $\max_{u_1\in[\beta_{\min},\beta_{\nom}]}\dot x_2(t) = \beta_{\nom} x_1(t) x_3(t) - \eta x_2(t) \leq \frac{\eta \hat x_2(\hat t)}{\hat x_3(\hat t)} x_3(\hat t) - \eta x_2(t) = \eta(\hat x_2(\hat t) - x_2(t)) \leq 0$ for all $t\in[\hat t_1, \hat t_2]$. 
	Thus, over $[\hat t_1, \hat t_2]$, $\hat x_2$ is decreasing, $\hat x_3$ is constant, and $\hat u_2$ is decreasing. 
	At $\hat t_2$, $\hat u_2(\hat t_2) = \gamma_{\nom}$, and we can again estimate this time using the bound on $\hat x_1$,
	\[
		\hat x_1(t) \leq e^{-\beta_{\min} \hat x_3(\hat t)(t - \hat t_1)}\hat x_1(\hat t_1),\quad t\in[\hat t_1, \hat t_2],
	\] 
	By exactly the same argument as before, $\hat t_2$ is finite. 
	
	At $\hat t_2$ we have $\dot{\hat x}_2(\hat t_2) + \dot{\hat x}_3(\hat t_2) = \beta_{\nom}\hat x_1(\hat t_2) \hat x_3(\hat t_2) - \gamma_{\nom}\hat x_3(\hat t_2)$, and $\dot{\hat x}_2(\hat t_2) + \dot{\hat x}_3(\hat t_2) \leq 0$. 
	Thus, $\beta_{\nom}\hat x_1(\hat t_2) \hat x_3(\hat t_2) - \gamma_{\nom}\hat x_3(\hat t_2) \leq 0$, which implies $\hat x_1(\hat t_2) \leq \frac{\gamma_{\nom}}{\beta_{\nom}} = \overline x_1$. 
	Moreover, because $\dot{\hat x}_3(\hat t_2) = \eta \hat x_2(\hat t_2) - \gamma_{\nom} \hat x_3(\hat t_2) = 0$ and $\hat x_3(\hat t_2) \leq I_{\max}$, $\hat x_2(\hat t_2) \leq \frac{\gamma_{\nom}I_{\max}}{\beta_{\nom}} = \overline x_2$. 
	Thus, $\x(\hat t_2)\in\XM$. 
	Finally, we use the control $\bold u_{\nom}$ for all $t\in[\hat t_2,\infty)$, which asymptotically drives the state to an equilibrium point. 
	
	Thus, we have constructed an input that drives the state into the invariant box $\XM$ in finite time, from any initial point in $\A^{\prime}$. 
	The integral of the continuous running cost, $\ell$, over the interval $[0,\hat t_2]$ is finite, and the cost $J_{\infty}(\x(\hat t_2))$ is finite, from Lemma~\ref{lem:J_inf_uni_bounded}. 
	Thus, there exists a $\bar C< \infty$ such that $V_{\infty}(\x^0) \leq \bar C$ for all $\x^0\in\A^{\prime}$, which completes the proof.  
%
%
%
%
%
%
%
%
%
\end{proof} 
We have collected enough facts regarding $\OCPT$ to now state the paper's main result. 
It says that if $T$ is sufficiently long, then the value function experiences a Lyapunov-like decrease on every iteration of the MPC loop. 
Thus, the state is driven to where $V_T(\x) = 0$ (i.e., a point on $\E$), and so $\A^{\prime}$ is rendered a domain of attraction of $\mathcal{E}_{\nom}$. 
\begin{theorem}\label{prop:main}
	For any control horizon, $\delta > 0$, there exists a finite prediction horizon, $T = N\delta$, $N\in\mathbb{N}$ with $N\geq 1$, such that the following relaxed Lyapunov inequality holds:
	\[
		V_T(\x(\delta;\x^0,\boldu^{\star}))\leq V_T(\x^0) - (1 - \alpha)\int_0^{\delta}\ell(\x(s;\x^0,\boldu^{\star}), \boldu^{\star}(s))\,\,\dee s,\,\,\forall \x^0\in\A^{\prime},
	\]
	where $\A^\prime$ is defined in \eqref{eq:def_set_A_prime} and $\alpha\in(0,1)$.
\end{theorem}
\begin{proof}
	We adapt the proof of \cite[Thm.~1]{esterhuizen2020recursive} to the setting of this paper, as we have a positive semi-definite running cost, $\ell$. 
	We highlight the parts of the proof that change, only sketching those parts that do not, and emphasise where the previous results established in this paper are used. 
	
	Recall from \eqref{quadratic_ell} that $\ell^\star(\x) = \lambda (x_2^2 + x_3^2)$. 
    With $\epsilon > 0$, let,
    \[
        \mathcal{N}_{\epsilon}(\E_{\nom}) := \{\x\in\mbbR^3 : \dist(\x, \E_{\nom}) < \epsilon \},
    \]
    where $\dist(\x, \E_{\nom}) := \inf_{\boldsymbol y \in \E_{\nom}} \Vert \x - \boldsymbol y \Vert$ is the distance function,
    be a small enough neighbourhood of the equilibrium points $\E_{\nom}$ such that the constant
	\[
	M:= \inf_{\x \in \X \setminus \mathcal{N}_{\epsilon}(\E_{\nom})} \ell^\star(\x) > 0,
	\]
	is well-defined. 
	By Proposition~\ref{lem:V_inf_uni_bounded_on_A} there exists a $C^\prime \geq 0$ such that $\A^\prime \subseteq V_{\infty}^{-1}[0,C^\prime] \subseteq V_{T}^{-1}[0,C^\prime]$ and we may conclude that,
	\[
		V_T(\x^0) \leq \frac{C^\prime}{M}\ell^\star(\x^0),\quad \forall \x^0 \in \A^\prime \setminus \mathcal{N}_{\epsilon}(\E_{\nom}).
	\]
	By Proposition~\ref{prop:cost_control},
	\[
	  V_{T}(\x^0)\leq V_{\infty}(\x^0) \leq \rho \ell^\star(\x^0),\quad \x^0 \in \mathcal{N}_{\epsilon}(\E_{\nom})\cap \XM.
	\]
	Combining these two inequalities we may deduce that,
	\[
		V_T(\x^0) \leq \sigma \ell^\star(\x^0),\quad \forall \x^0 \in \A^\prime,
	\]
	where $\sigma := \max\{\frac{C^\prime}{M}, \rho\}$. 
	By the same arguments as in the proof of \cite[Thm.~1]{esterhuizen2020recursive} we get, for any $\x^0 \in \A^\prime$,
	\begin{equation}
		\min\{\sigma \ell^\star(\x^0), C^\prime\} \geq V_T(\x^0) \geq \left(\frac{\sigma + \delta}{\sigma}\right)^{N - 1} \int_{T - \delta}^T \ell(\x(s;\x^0, \boldu^\star), \boldu^\star)\, \dee s.\label{ineq_l_M}
	\end{equation}
	If $N$ is chosen large enough such that
	\begin{equation}
		\frac{C^\prime}{M \delta}\left(\frac{\sigma}{\sigma  + \delta}\right)^{N - 1} < 1,\label{N_big_1}
	\end{equation}
	then, referring to \eqref{ineq_l_M},
	\[
		\int_{T - \delta}^T \ell(\x(s;\x^0, \boldu^\star), \boldu^\star)\, \dee s < M \delta.
	\]
	By exactly the same arguments as in \cite[Thm.~1]{esterhuizen2020recursive} it can then be shown that, for any $\x^0 \in \A^\prime$,
	\begin{equation}
		V_T(\x(\delta;\x^0,\boldu^{\star}))\leq V_T(\x^0) - \int_0^\delta \ell(\x(s;\x^0, \boldu^\star), \boldu^\star)\, \dee s + \frac{\sigma^2}{\delta}\left(\frac{\sigma}{\sigma  + \delta}\right)^{N - 1} \ell^\star(\x^0). \label{eq:decrease_V}
	\end{equation}
	We now argue that there exists a $\bar C \geq 0$ such that,
	\[
		\bar C \int_0^\delta \ell(\x(s;\x^0, \boldu^\star), \boldu^\star)\, \dee s \geq \delta \ell^\star(\x^0),\quad \forall \x^0 \in\A^\prime.
	\]
	Recall from the dynamics that,
	\begin{align}
		\dot x_2(t) \; & = \; \beta(t) x_1(t) x_3(t) - \eta x_2(t) \; \geq \; - \eta x_2(t), \notag \\
		\dot x_3(t) \; & = \; \eta x_2(t) - \gamma(t) x_3(t) \; \geq \; - \gamma_\mathrm{max} x_3(t). \notag 
	\end{align}
	Thus, from Grönwall's lemma,
	\begin{align*}
		\int_0^\delta \ell(\x(s;\x^0, \boldu^\star), \boldu^\star)\, \dee s \; & = \; \min_{\boldu\in\U_{\delta}(\x^0)} \int_0^{\delta}\Vert \x(s;\x^0,\boldu) \Vert_Q^2 + \Vert \boldu(s) - \boldu_{\nom}\Vert_R^2 \,\, \dee s \; \\
		& \geq \; \int_0^{\delta} e^{-2\eta s}(x_2^0)^2 + e^{-2 \gamma_\mathrm{max} s}(x_3^0)^2 \, \dee s \\
		& \geq \; \mathrm{e}^{-2 \eta \delta} \int_0^\delta (x_2^0)^2 + \mathrm{e}^{-2 \gamma_\mathrm{max} \delta} \int_0^\delta (x_3^0)^2 \, \dee s \\
		& \geq \; \mathrm{e}^{- 2 \delta \max \{ \eta, \gamma_\mathrm{max} \}} \delta \ell^\star(\x^0).
	\end{align*}	
	Or, equivalently, 
	\begin{align}
		\bar C \int_0^\delta \ell(\x(s;\x^0, \boldu^\star), \boldu^\star)\, \dee s \; \geq \; \delta \ell^\star(\x^0),\quad \forall\x^0\in\A^{\prime}, \nonumber 
	\end{align}
	with $\bar C := \max_{\delta \in (0,T]} \mathrm{e}^{2 \delta \max \{ \eta, \gamma_\mathrm{max} \}} = \mathrm{e}^{ 2 T \max \{ \eta, \gamma_\mathrm{max} \}}$.
	Referring back to \eqref{eq:decrease_V}, we get, for $\x^0 \in\A^\prime$,
	\[
		V_T(\x(\delta;\x^0,\boldu^{\star}))\leq V_T(\x^0) - 
		(1- \alpha) \int_0^\delta \ell(\x(s;\x^0, \boldu^\star), \boldu^\star)\, \dee s,
	\]
	where
	\[
		\alpha := \bar C\frac{\sigma^2}{\delta^2}\left( \frac{\sigma}{\sigma + \delta} 	\right)^{N - 1} \geq 0.
	\]
	If $N$ is now chosen large enough such that, $\alpha < 1$, then $V_T$ experiences a decrease at $t = \delta$. 
	Thus, referring back to \eqref{N_big_1}, we have shown that the statement in the proposition holds if $N$ is chosen such that
	\[
		\max \biggl\{ \frac{C^\prime}{M \delta}\left(\frac{\sigma}{\sigma  + \delta}\right)^{N - 1}, \underbrace{\bar C\frac{\sigma^2}{\delta^2}\left( \frac{\sigma}{\sigma + \delta} 	\right)^{N - 1}}_{\alpha} \biggr\} < 1. 
	\]
	This completes the proof. 
\end{proof}

\section{Numerical simulation}\label{sec:numerics}

We consider the constrained SEIR system, $\mathcal{S}$, where $f$ is specified in \eqref{eq:SIER_f}, with the following parameters, $\beta_{\nom} = 0.44$, $\gamma_{\nom} = 1/6.5$, $\beta_{\min} = 0.22$, $\gamma_{\max} = 0.5$, $\eta = 1/4.6$ and $I_{\max} = 0.05$. 
These are the same constants considered in the numerics of \cite{sauerteig2022model}. 
We solve $\OCPT$ on every iteration of the loop in Algorithm~\ref{alg:mpc} via a direct approach. That is, we discretise the dynamics via the explicit Euler scheme, see for example \cite{Gerdts2012}, with a step size of $h = 0.25$ days and solve the associated nonlinear programme with the interior-point method of Matlab's $\mathsf{fmincon}$ function. 
By using the Euler scheme the solution of the optimisation problem could be sped up considerably. 
Moreover, the error in integration of the dynamics between the Euler scheme with $h=0.25$ and the Runge-Kutta 4 scheme with very small constant step size was acceptably small. 

\begin{figure}[htb]
	\centering
	\includegraphics[scale=0.6]{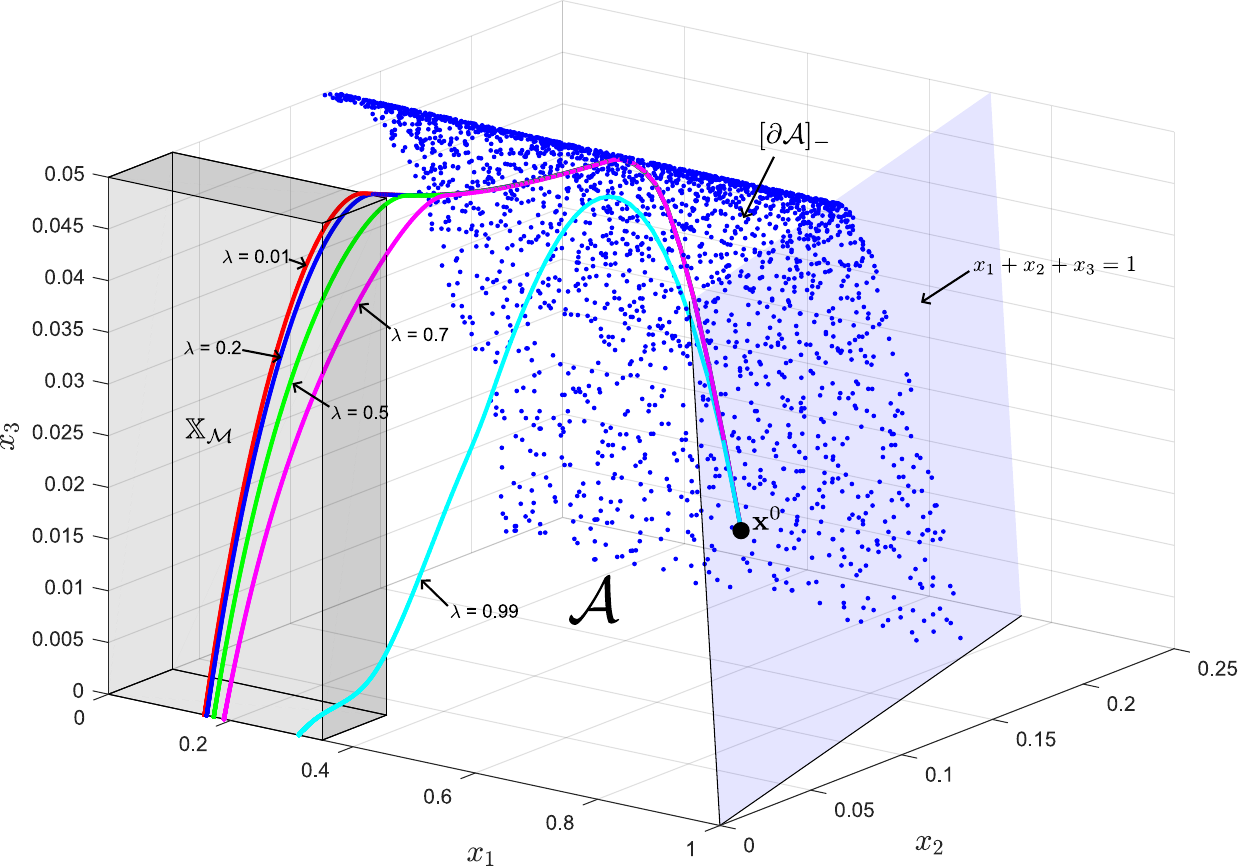}
	\caption{Closed-loop state trajectories of the constrained SEIR system, $\mathcal{S}$, produced by the MPC algorithm without terminal ingredients, Algorithm~\ref{alg:mpc}, with varying values for $\lambda$ (0.99 cyan, 0.7 magenta, 0.5 green, 0.2 dark blue, 0.01 red). 
		The blue points indicate the boundary of the admissible set,  $[\partial \A]_- := \partial \A \cap \mathrm{int}(\X)$; the grey box indicates the invariant box, $\XM$; and the blue plane indicates the set $\{\x\in\mbbR^3: x_1 + x_2 + x_3 = 1\}$. }
	\label{fig:3D_plot}
\end{figure}

Figure~\ref{fig:3D_plot} shows the closed loop trajectories produced via Algorithm~\ref{alg:mpc} with initial state $\x^0 = (0.5, 0.18, 0.01)^\top\in\A^{\prime}$, control horizon $\delta = 1$ day, prediction horizon $T=20$ days, and $\lambda = 0.99$ (cyan), $\lambda = 0.7$ (magenta), $\lambda = 0.5$ (green), $\lambda = 0.2$ (dark blue) and $\lambda = 0.01$ (red). 
To solve $\OCPT$ with these parameters requires roughly 10 seconds on a standard desktop computer. 
Thus, for a simulation of the closed loop over 100 days (a window over which the most important aspects of the closed loop behaviour is observed), this equates to about 16 minutes of computation time. 
The blue dots, produced via the analysis in \cite{ester_epidemic_management_2021}, are on the so-called \emph{barrier}, labelled $[\partial \A]_- := \partial \A \cap \mathrm{int}(\X)$. 
This is the part of the boundary of the admissible set that is also in the interior of $\X$. 
The grey box indicates the box $\XM$, and the blue plane indicates the set $\{\x\in\mbbR^3: x_1 + x_2 + x_3 = 1\}$. 
We terminate the simulation once $ \max\{x_2(t),x_3(t) \} \leq 10^{-8}$ (in a population of 100 million people, this corresponds to no one being infected).

\begin{figure}[htb!]
	\centering
	\includegraphics[scale=0.75]{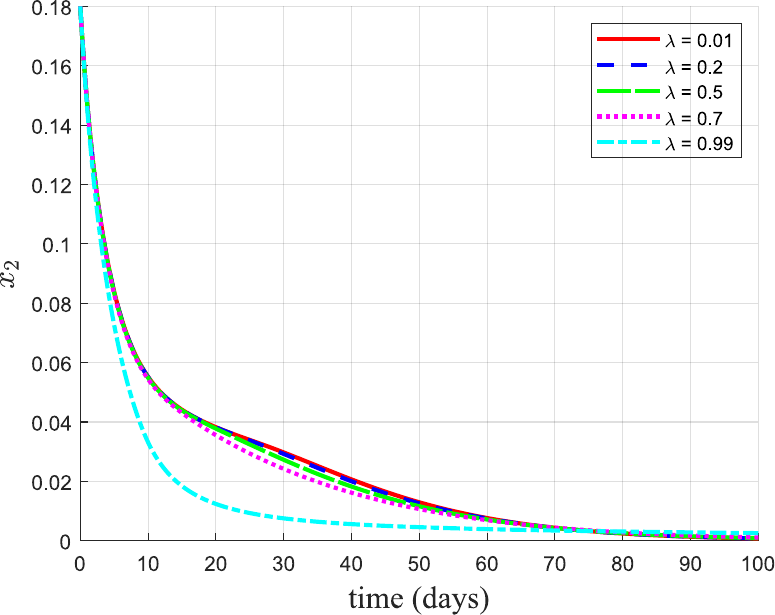}
	\caption{Closed-loop behaviour of the exposed compartment under the MPC feedback with varying~$\lambda$ (first 100 days).}
	\label{fig:E_plot}
\end{figure}

\begin{figure}[htb!]
	\centering
	\includegraphics[scale=0.75]{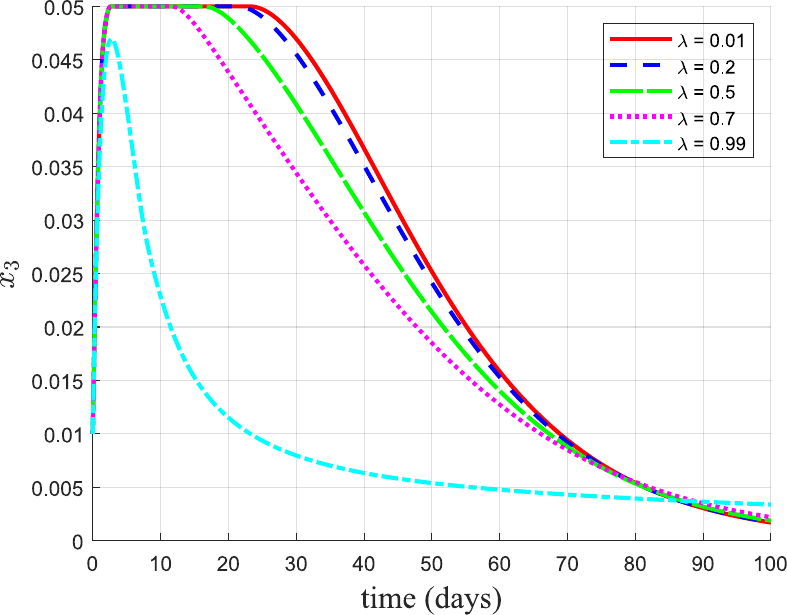}
	\caption{Closed-loop behaviour of the infectuous compartment under the MPC feedback with varying~$\lambda$ (first 100 days).}
	\label{fig:I_plot}
\end{figure}

\begin{figure}[htb!]
	\centering
	\includegraphics[scale=0.75]{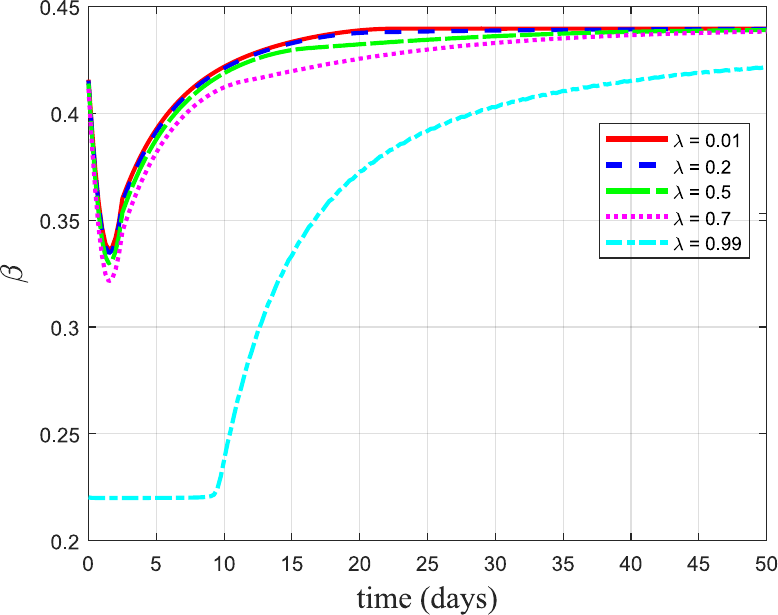}
	\caption{Closed-loop behaviour of $\beta$ with varying~$\lambda$ (first 50 days).}
	\label{fig:beta_plot}
\end{figure}

\begin{figure}[htb!]
	\centering
	\includegraphics[scale=0.75]{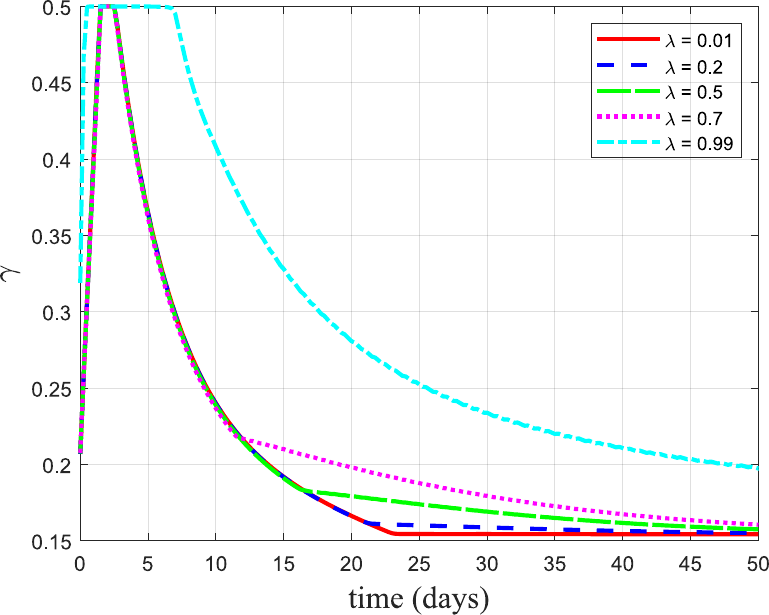}
	\caption{Closed-loop behaviour of $\gamma$ with varying~$\lambda$ (first 50 days).}
	\label{fig:gamma_plot}
\end{figure}

Figures~\ref{fig:E_plot} and \ref{fig:I_plot} show the closed-loop exposed and infected compartments, respectively, as a function of time, whereas Figures~\ref{fig:beta_plot} and \ref{fig:gamma_plot} show the closed-loop control inputs, $\beta$ and $\gamma$, as a function of time.
As clearly seen in Figures~\ref{fig:3D_plot}-\ref{fig:gamma_plot}, larger values for $\lambda$ result in fewer total infectuous and exposed individuals during the epidemic, but with the requirement of stricter societal interventions. 
However, perhaps counter to intuition, larger values for $\lambda$ (which penalises deviation of the state from equilibrium points) also result in the pandemic lasting longer, see Table~\ref{table:pandemic_lifespan}. 

It is worth noting that, with $\lambda = \frac{1}{2}$, the shortest horizon for which $\OCPT$ was recursively feasible and for which the closed loop was driven to an equilibrium point was $T = 2$ days. 
Compare this to the longer horizon of $T = 25$ days that was required in the MPC scheme \emph{with} terminal ingredients in \cite{sauerteig2022model}. 
The reason for this is that the value function of $\OCPT$ (i.e. the problem without terminal ingredients) experiences the Lyapunov-like decrease in Theorem~\ref{prop:main} well before the state reaches the invariant box, $\XM$, which is the target set in the MPC scheme from \cite{sauerteig2022model}. 
\begin{table}[htb!]
	\caption{The effect of $\lambda$ on the lifetime of the epidemic.}\label{table:pandemic_lifespan}
	\centering
	\begin{tabular}{|c|cccc|}
		\hline
		& \multicolumn{4}{c|}{Days until $\max\{x_2(t), x_3(t) \} < \mathrm{val}$} \\
		\hline
		$\lambda$ & $\mathrm{val}=10^{-5}$ & $\mathrm{val}=10^{-6}$ & $\mathrm{val}=10^{-7}$ & $\mathrm{val}=10^{-8}$\\
		\hline
		0.01 & 186.5 & 225 & 263.75 & 302 \\
		\hline
		0.2 &  188.75 & 228 & 267.5 & 306.5 \\
		\hline
		0.5 & 196.75 & 239 & 281.25 & 323.75 \\
		\hline
		0.7 & 212.25 & 260 & 307.5 & 355.25 \\
		\hline
		0.99 & 612 & 801.25 & 988 & 1175\\
		\hline
	\end{tabular}
\end{table}

\section{Conclusions and Perspectives}\label{sec:conclusion}

Drawing on results from \cite{esterhuizen2020recursive}, we rigorously showed that an SEIR epidemic can be eliminated via MPC without terminal ingredients provided that the state initiates from a suitable subset of the admissible set, one uses a suitable quadratic running cost, and the prediction horizon in the finite-horizon OCP is sufficiently long. 
Numerics suggest that this prediction horizon may be much shorter than the one required in MPC with terminal ingredients. 
Moreover, they showed that higher values of $\lambda$ in the cost functional, which penalise deviation from equilibrium points, result in fewer infected individuals but longer epidemics. 

Future research could generalise the results to the case where the running cost, \eqref{quadratic_ell}, is more elaborate, allowing certain state variables or intervention measures to carry more weight than others. 
We note that the use of a general quadratic cost (that is, one of the form $\ell(\x,\boldu) = \Vert \x \Vert_Q^2 + \Vert \boldu -\boldu_{\nom}\Vert_R^2$ with $Q\succeq 0$ and $R\succ 0$) introduces cross terms, which make the proofs of some of the results more difficult. 

Future research could also focus on generalisations of the paper's results to more elaborate epidemic models. 
In particular, models where \emph{reinfection} is possible (modelled, for example, as a flow from the removed compartment, $R$, to the susceptible compartment, $S$) are interesting for diseases like Covid-19. 
The main challenge introduced by such models is the fact that the susceptible compartment is not strictly decreasing over time, a fact we used extensively in deriving the results in this paper. 



%
%
%

\appendix

\section{Proof of Lemma~\ref{lem:J_inf_uni_bounded}}

\begin{proof}
	\sloppy	
	To ease notation throughout the proof, let $\x(t) = (x_1(t), x_2(t), x_3(t))^\top = \x(t;\x^0,\boldu_{\nom})$ denote the ``nominal'' solution, and let $x_i^{\infty} = x_i^{\infty}(\x^0,\boldu_{\nom})$, $i\in{1,2,3}$. 
	Because $\x^0\in\XM$, by Lemma~\ref{lem:facts_of_sets}, we have that $\x(t)\in\X$ for all $t\geq 0$, with $x_1^{\infty}\in [0,1]$, and $x_2^{\infty} = x_3^{\infty} = 0$. 
	Therefore, noting that,
	\[
	\int_0^{\infty} \beta_{\nom}x_1(t)x_3(t)\,  \dee t = -\int_0^{\infty} \dot{x}_1(t)\,  \dee t = x_1^0- x_1^{\infty}\in[0,1],
	\]
	we see that,
	\begin{align*}
		\int_0^{\infty}x_2(t)\,  \dee t & = \frac{1}{\eta}\left(\int_0^{\infty} \beta_{\nom} x_1(t) x_3(t)\,  \dee t- \int_0^{\infty}\dot{x}_2(t)\right)\,  \dee t \\
		& = \frac{1}{\eta}(x_1^0 - x_1^{\infty} + x_2^0) < \infty.
	\end{align*}
	Similarly,
	\begin{align*}
		\int_0^{\infty}x_3(t)\,  \dee t & = \frac{1}{\gamma_{\nom}}\left(\int_0^{\infty} \eta x_2(t) \,  \dee t - \int_0^{\infty} \dot{x}_3(t)\,  \dee t \right)\\
		& = \frac{1}{\gamma_{\nom}}\left( \eta \frac{1}{\eta}(x_1^0 - x_1^{\infty} + x_2^0) + x_3^0 \right) < \infty.
	\end{align*}	
	Let $C := \lambda\left( \frac{1}{\eta}(x_1^0 - x_1^{\infty} + x_2^0) + \frac{1}{\gamma_{\nom}}(x_1^0 - x_1^{\infty} + x_2^0 + x_3^0 )\right)<\infty$. 
	Then, because $x_2(t),x_3(t)\in[0,1]$ for all $t\geq 0$, we see that,
	\begin{align*}
		J_{\infty}(\x^0,\boldu_{\nom}) & = \lambda \int_0^{\infty} \left(x_2^2(t) + x_3^2(t)\right)\,  \dee t\\
		& \leq \lambda \int_0^{\infty} \left( x_2(t) + x_3(t) \right)\  \dee t \\
		& \leq C,\quad \forall \x^0\in\XM.
	\end{align*}
\end{proof}

  \bibliographystyle{elsarticle-num} 
  \bibliography{references_MPC_without_stab_epidemic_models}


%
%
%

\end{document}